ANNALES
DE L'INSTITUT
HENRI
POINCARÉ
PROBABILITÉS
ET STATISTIQUES

# A lower bound for the principal eigenvalue of the Stokes operator in a random domain[*]

## V. V. Yurinsky

*Departamento de Matemática, Universidade da Beira Interior, Rua Marquês d'Ávila e Bolama, 6201-001 Covilhã, Portugal.*
*E-mail: yurinsky@ubi.pt*

**Abstract.** This article is dedicated to localization of the principal eigenvalue (PE) of the Stokes operator acting on solenoidal vector fields that vanish outside a large random domain modeling the pore space in a cubic block of porous material with disordered micro-structure. Its main result is an asymptotically deterministic lower bound for the PE of the sum of a low compressibility approximation to the Stokes operator and a small scaled random potential term, which is applied to produce a similar bound for the Stokes PE. The arguments are based on the method proposed by F. Merkl and M. V. Wütrich for localization of the PE of the Schrödinger operator in a similar setting. Some additional work is needed to circumvent the complications arising from the restriction to divergence-free vector fields of the class of test functions in the variational characterization of the Stokes PE.

**Résumé.** Cet article est dédié à l'étude de la localisation de la valeur propre principale (VPP) de l'opérateur de Stokes sous la condition de Dirichlet sur la frontière d'un grand domaine aléatoire qui modélise l'espace des pores d' un bloc cubique de matière poreuse dotée d'une microstructure désordonnée. Le résultat principal est une borne inférieure asymptotiquement déterministe pour la VPP de l'opérateur correspondant á l'écoulement d'un liquide peu compressible en présence d'un petit potentiel positif aléatoire. Les arguments sont fondés sur la méthode proposée par F. Merkl et M. V. Wütrich pour localiser la VPP de l'opérateur de Schrödinger dans une situation similaire. Des efforts supplémentaires sont nécessaires pour combattre les complications provenant de la réduction à la classe de champs vectoriels de divergence nulle de la famille des fonctions utilisées pour caractériser la VPP de l'opérateur de Stokes par une formule variationnelle.

[*]Work supported by FCT (Portugal) through Centro de Matemática UBI, Project DECONT – *Differential Equations of Continuum Mechanics* in the framework of program POCI 2010 co-financed by the Portuguese Government and EU (FEDER).

## 1. Introduction

This article deals with localization of the principal eigenvalue (PE) of the Stokes operator acting on solenoidal vector fields over a fine-grained random domain that models the pore space in a large block of a material with disordered micro-structure (e.g., porous rock). Below, the flow domain is

$$F_t = \widehat{Q}_t^0 \setminus S, \quad \widehat{Q}_t^0 = \left(-\frac{1}{2}t, \frac{1}{2}t\right)^d \subset \mathbb{R}^d, \quad t \gg 1, \tag{1.1}$$

where the closed set $S = S(\omega) = \{x: \ V(x, \omega) = 1\} \subseteq \mathbb{R}^d$ models the "skeleton" of the porous material, and $V$ is a measurable $\{0, 1\}$-valued random field.

The PE $\mathfrak{S}_t$ of the Stokes operator acting on solenoidal velocity fields from the Sobolev space $\mathbf{H}_0^1(F_t)$ admits the following well-known variational characterization in terms of Rayleigh quotients (see [8], Chapter 1.8):

$$\mathfrak{S}_t = \inf\{\|\phi\|_2^{-2}\|\nabla\phi\|_2^2 \colon \ \mathrm{div}(\phi) = 0, \phi \in \mathbf{H}_0^1(F_t)\}. \tag{1.2}$$

Under natural assumptions about the structure of the flow domain, this PE exhibits essentially deterministic asymptotic behaviour as $t \to \infty$.

Rigorous results on asymptotically deterministic behaviour of the PE of an elliptic operator with random elements originate in the work of A.-S. Sznitman on localization of the PE of the Laplacian under the Dirichlet condition on the boundary of a random domain (see [6, 7] and the bibliography therein).

The method of enlargement of obstacles [7], used in most earlier publications to derive a lower bound on the PE, compares it with the Dirichlet PE's for subdomains of simpler shape that are compatible with a typical configuration of the random element (see, e.g., [6, 7, 9]).

Later, Merkl and Wütrich [4] elaborated a new method to localize the PE of the Schrödinger operator with a "scaled" small random non-negative potential term by analyzing feasibility of specific values of the Rayleigh quotient for individual test functions. This article adapts the approach of [4] to flows in porous media.

*Lower bound on Stokes PE*

To show that the large-volume asymptotic behavior of PE (1.2) is essentially deterministic, it suffices to find for it asymptotically equivalent deterministic upper and lower confidence bounds. This can be done by techniques quite similar to those used for localization of the Laplacian's PE in the same setting.

Yet, the use of divergence-free fields as test functions in (1.2) complicates the construction of a confidence interval for the Stokes PE and makes it less explicit – both unilateral bounds include the constant

$$\mathcal{S} = \inf\{\|\phi\|_2^{-2}\|\nabla\phi\|_2^2 \colon \ \mathrm{div}(\phi) = 0, |\{|\phi| > 0\}| = 1, \phi \in \mathbf{H}^1(\mathbb{R}^d)\}, \tag{1.3}$$

which[1] can be loosely interpreted as the smallest value that the PE of the Stokes operator can have for domains of unit measure. It is strictly positive and at least as large as the Faber–Krahn bound for the Dirichlet PE of the Laplacian because of the additional restriction on the class of test functions. Constant (1.3) can be approximated ([10], Lemma A.1) by its counterparts for the low compressibility approximations to the Stokes operator [8], Chapter 1.6:

$$\mathcal{S} = \lim_{\alpha \to 0+} \mathcal{C}_\alpha, \quad \mathcal{C}_\alpha = \inf\{\|\phi\|_2^{-2}\mathcal{K}_\alpha(\phi) \colon \ |\{|\phi| > 0\}| = 1\}, \tag{1.4}$$

where $\mathcal{K}_\alpha(\phi) = \int_{\mathbb{R}^d} \widehat{K}_\alpha(\phi(x))\,\mathrm{d}x$, $\widehat{K}_\alpha(\phi) = |\nabla\phi|^2 + \alpha^{-1}(\mathrm{div}(\phi))^2$, and the vector-valued test functions are from $\mathbf{H}^1(\mathbb{R}^d)$.

Only the lower confidence bound for PE (1.2) is derived below – the approach of [4] is applied to prove the following theorem of [10].

---

[1]Here and below $|A|$ is the Lebesgue measure of a set $A \subset \mathbb{R}^d$. Constant (1.3) was introduced in [10], Eq. (1.10), in a different form.

**Theorem 1.1.** *Denote by $S_z(\omega) = S \cap \{x : x - z \in (-\frac{1}{2}, \frac{1}{2}]^d\}$ the complement to the flow region in the unit cube centered at $z \in \mathbb{Z}^d$.*

*If there exist independent identically distributed random variables $\xi_z$, $z \in \mathbb{Z}^d$, such that*

$$0 \leq \xi_z \leq 1, |S_z| \geq \xi_z, \quad \mathbf{P}\{\xi_z = 0\} \leq p, \qquad \mu = \mathbf{E}\xi_z > 0, \tag{1.5}$$

*then for each $\varepsilon > 0$*

$$\lim_{t \to \infty} \mathbf{P}\left\{(\ln t)^{2/d}\mathfrak{S}_t > \left(\frac{\nu}{d}\right)^{2/d}\mathcal{S} - \varepsilon\right\} = 1, \quad \nu = \ln\left(\frac{1}{p}\right). \tag{1.6}$$

The matching upper bound

$$\forall \varepsilon > 0 \quad \lim_{t \to \infty} \mathbf{P}\left\{(\ln t)^{2/d}\mathfrak{S}_t < \left(\frac{1}{d}\ln\left(\frac{1}{p}\right)\right)^{2/d}\mathcal{S} + \varepsilon\right\} = 1 \tag{1.7}$$

was derived in [10] (for $d = 2$) and [11] (for $d \geq 3$) for a model of porosity where the skeleton consists of isolated components. In this model, the indicator function of the set $S$ in (1.1) satisfies the inequality

$$1_S(x; \omega) \leq \sum_{z \in \mathbb{Z}^d} \varepsilon_z(\omega) 1_W(x - z), \tag{1.8}$$

the closed set $W \Subset (-\frac{1}{2}, \frac{1}{2})^d$, $|W| > 0$, is sufficiently regular, and the binary random variables $\varepsilon_z \in \{0, 1\}$ are independent and identically distributed, $p = \mathbf{P}\{\varepsilon_z = 0\} = 1 - \mathbf{P}\{\varepsilon_z = 1\} \in (0, 1)$. Combined, Eqs (1.6) and (1.7) show that conditions (1.5) and (1.8) ensure deterministic asymptotic behaviour of $\mathfrak{S}_t$ for this model of porosity, and it is not affected by translations and rotations of the skeleton.

The bounds (1.6) and (1.7) are, unfortunately, much less explicit then the well-known results on localization of the PE of the Laplacian [7] using the Faber–Krahn inequality.

It is obvious that in definition (1.3) a minimizer, if it exists, cannot be unique because the problem is invariant with respect to translations and rotations. Both calculation of $\mathcal{S}$ and characterization of shapes of the sets whose Stokes PE's are close to this constant seem to be open problems. The proof of (1.4) mentioned above does not provide any practical approach to calculating $\mathcal{S}$.

In the proof of Theorem 1.1, information about the possible shapes of sets with small Stokes PE's is relatively unimportant.

By contrast, the main difficulty in the proof of (1.7) lies in showing that a typical configuration of the skeleton allows the existence of divergence-free test functions with bounded support and Rayleigh quotients sufficiently close to $\mathcal{S}$.

For $d = 2$, the choice of test functions in definition (1.3) can be limited to ones having bounded support ([10], Lemma A.4) whatever the configuration of the skeleton. For $d \geq 3$, this question remains open, and the derivation of the upper bound in [11], Lemma 2.1, is more complicated because it exploits the possibility of adjusting a divergence-free test function with low Rayleigh quotient to a given configuration of the flow region whenever this contains a sufficiently large "vacuity" (a connected subset free from inclusions of the skeleton).

*Low compressibility approximation to Stokes PE*

The restriction of the class of admissible test functions in (1.2) to divergence-free ones precludes the use of some techniques of [4] that employ cutoffs.

To bypass this difficulty, the method of [4] is used to derive the lower bound first for the Dirichlet PE of the auxiliary operator that acts on smooth functions as

$$\Lambda_{\alpha,\beta,t}\phi = -\left(\Delta\phi + \left(\frac{1}{\alpha}\right)\nabla\operatorname{div}(\phi)\right) + (\ln t)^{-2/d}\beta V\phi. \tag{1.9}$$

It combines a low-compressibility approximation to the Stokes operator ([8], Chapter 1.6) with a small potential term that substitutes the random boundary [4]. For a given configuration of skeleton, the Dirichlet PE of operator (1.9) is

$$\mathfrak{C}_{\alpha,\beta,t} \overset{\text{def}}{=} \inf_{\phi \in \mathbf{H}_0^1(\widehat{Q}_t^0)} \|\phi\|_2^{-2}(\mathcal{K}_\alpha(\phi) + (\ln t)^{-2/d}\beta \|V^{1/2}\phi\|_2^2), \tag{1.10}$$

where notation is that of (1.4) and the test functions are extended from $\widehat{Q}_t^0$ to all $\mathbb{R}^d$ by zero. It is obvious from definitions (1.2) and (1.10) that

$$\mathfrak{S}_t \geq \mathfrak{C}_{\alpha,\beta,t}. \tag{1.11}$$

Technically, the main result of this article is the following theorem on the limit behaviour of PE (1.10) under normalization

$$\lambda(t,\alpha,\beta) \overset{\text{def}}{=} (\ln t)^{2/d}\mathfrak{C}_{\alpha,\beta,t}. \tag{1.12}$$

**Theorem 1.2.** *If condition (1.5) is satisfied, then the normalized PE (1.12) admits the deterministic confidence bounds*

$$\forall \varepsilon > 0 \quad \lim_{t \to \infty} \mathbf{P}\{\lambda(t,\alpha,\beta) > \mathcal{C}_{\alpha,\beta} - \varepsilon\} = 1, \tag{1.13}$$

*where in notation of (1.10)*

$$\mathcal{C}_{\alpha,\beta} \overset{\text{def}}{=} \inf\{\mathcal{K}_\alpha(\phi) + \beta\mathcal{G}(\phi;d): \ \phi \in \mathbf{H}^1(\mathbb{R}^d), \|\phi\|_2 = 1\},$$

*and the function $G(u) = -\ln \mathbf{E}\exp\{-u\xi_0\}$ is used to define*

$$\Gamma(h;\phi) = \int_{\mathbb{R}^d} G(h|\phi(x)|^2)\,\mathrm{d}x, \quad h > 0, \tag{1.14}$$

$$\mathcal{G}(\phi;D) = \sup_{h>0}\{h^{-1}(\Gamma(h;\phi) - D)\}, \quad D \in \mathbb{R}. \tag{1.15}$$

Theorem 1.1 is deduced from Theorem 1.2 and inequality (1.11).

Theorem 1.2 and its proof mainly digress from [4] in the form of the functional characterizing the feasible values of individual Rayleigh quotients and the use of Cramér's transform to derive a tractable estimate for the exponential moment of $\|V^{1/2}\phi\|_2^2$.

The proofs of Theorem 1.2 and Theorem 1.1 occupy, respectively, Sections 2 and 3. The appendices contain some necessary auxiliary material.

*Notation*

Points in $\mathbb{R}^d$ and their coordinates are denoted $x = (x_j)$. The scalar product is $x \cdot y = \sum x_j y_j$ and $|x| = \sqrt{x \cdot x}$ is the corresponding norm. One more norm in use is $|x|_* = \max_j |x_j|$, and the corresponding distance from a point to a set is $\mathrm{Dist}(x,B) = \inf\{\max_j |x_j - y_j|: y \in B\}$. For a $d \times d$ matrix $a = (a_{jk})$ the norm $|a| = (a: a)^{1/2}$ corresponds to the product $a: b = \sum_{j,k=1}^d a_{jk}b_{jk}$.

For a finite set, $\#(S)$ is the number of its elements, and $|A|$ is the Lebesgue measure of a set $A \subset \mathbb{R}^d$. Sets obtained by translations and changes of scale are denoted

$$a + \alpha G = \{x \in \mathbb{R}^d: \ x = a + \alpha y, y \in G\}, \quad a \in \mathbb{R}^d, \alpha \in \mathbb{R}, G \subseteq \mathbb{R}^d.$$

The cube $(-\frac{1}{2}, \frac{1}{2}]^d$ is always denoted $Q$; $Q^0$, $\overline{Q}$ are its interior and closure.

For a real-valued function $\nabla \psi = (\nabla_j \psi)$ is the gradient and $|\nabla \psi|$ its norm; by analogy $|\nabla \phi|^2 = \sum_{j,k=1}^{d} (\nabla_j \phi^{(k)})^2$ for a function $\phi = (\phi^{(j)}) \in \mathbb{R}^d$. The divergence is div($\phi$). Notation of integrals is often abbreviated: $\int_G f(x) \, dx$ may be reduced to $\int_G f$ or $\int f$ if the context excludes misunderstanding. $\mathbf{P}$ and $\mathbf{E}$ denote probability and expectation on the probability space $\langle \Omega, \mathcal{F}, \mathbf{P} \rangle$.

Notation of function spaces follows [8] and [1]. For $G \subseteq \mathbb{R}^d$, the spaces of scalar or vector valued summable functions are denoted $L^p(G)$, and $\| \cdot \|_p$ or $\| \cdot \|_{L^p(G)}$ is the corresponding norm. For a bounded open set $G$, the Sobolev space $H_0^1(G)$ is the closure in the norm $\|\phi\|_{H^1} = (\|\phi\|_2^2 + \|\nabla\phi\|_2^2)^{1/2}$ of the space $C_0^\infty(G)$ of scalar smooth functions with compact support in $G$. Its counterpart for vector valued functions is $\mathbf{H}_0^1(G)$, and as in [8], Chapter 1.1.4, the subspace of solenoidal fields on $G$ is $\mathbf{V}(G) = \{\phi \in \mathbf{H}_0^1(G): \text{div}(\phi) = 0\}$. The spaces $H^1(\mathbb{R}^d)$ and $\mathbf{H}^1(\mathbb{R}^d)$ are the closures of the set of compact-supported smooth scalar and vector valued functions in the norm $\|\phi\|_{H^1}$.

Positive constants are denoted $c, c_i, \hat{c}$, etc. No attempt is made to keep track of their numerical values, so the same notation may be used for different quantities depending on the context. Implicit "equalities" similar to $c \cdot c = c$, $c + c = c$, etc. mean that the value of a new constant appearing in a calculation is determined by the same parameters as those of the old ones.

The calculations below use some multiplicative inequalities for the Sobolev space $H_0^1(G)$ ([3], Chapter II.2):

$$\begin{aligned}
&\|\phi\|_q \leq C(q) \|\phi\|_2^{2/q} \|\nabla\phi\|_2^{1-2/q}, \quad q > 2, d = 2, \\
&\|\phi\|_{2/(1-2/d)} \leq C(d) \|\nabla\phi\|_2, \qquad d \geq 3.
\end{aligned} \tag{1.16}$$

One more tool is a modification of the Poincaré–Friedrichs inequality: if $\phi$ is an $H^1$ function defined on a convex set $Q^+$ and $Q^0, Q^1 \subseteq Q^+$ are its subsets such that $|Q^1|/|Q^0| \leq \alpha_*^d$, then

$$\|\phi\|_{L^2(Q^1)}^2 \leq 2\alpha_*^d \|\phi\|_{L^2(Q_0)}^2 + \rho^2(Q^+) C(\alpha_*) \|\nabla\phi\|_{L^2(Q^+)}^2, \tag{1.17}$$

where $\rho(Q^+)$ is the diameter of $Q^+$; the constant on the right-hand side admits the estimate $C(\alpha_*) \leq c(d) \max\{\alpha_*, \alpha_*^{d-1}\}$ with $c(d)$ that depends on the dimension $d$ alone (see [9] for a proof).

## 2. Low compressibility bound on PE

### 2.1. Reduction to smaller boxes

Below the original spatial variable $\hat{x} \in \hat{Q}_t$ is changed to

$$x = \tau \hat{x} \in Q_t \overset{\text{def}}{=} \tau \hat{Q}_t = \tau t Q^0, \quad \tau = (\ln t)^{-1/d}. \tag{2.1}$$

In the new variables, the eigenvalue problem for operator (1.9) becomes $-(\Delta\phi + \alpha^{-1}\nabla \text{div}(\phi)) + \beta V_t \phi = \lambda \phi$, $\phi|_{\partial Q_t} = 0$, where $V_t(x) = V(\tau^{-1}x; \omega)$. Its PE equals the normalized PE of (1.12) and (1.13):

$$\lambda(t, \alpha, \beta) = \inf\{\|\phi\|_2^{-2}(\mathcal{K}_\alpha(\phi) + \beta\|V_t^{1/2}\phi\|_2^2): \ \phi \in \mathbf{H}_0^1(Q_t)\}. \tag{2.2}$$

The parameter $\tau$ of (2.1), $\mu$ of (1.5), and a large odd integer number $T = T(t)$ are used below to define partitions of $\mathbb{R}^d$ into cells $C_z^0 = \tau(z + Q)$ of size $H_0 = \tau$ and blocks $C_\zeta^1 = H_1(\zeta + Q)$ with $H_1(t) = \tau T$. The sets

$$\begin{aligned}
&\mathbb{Q}_t = \{z \in \mathbb{Z}^d: \ C_1^1 \cap Q_t \neq \emptyset\}, \quad \#(\mathbb{Q}_t) \leq \left(\frac{t}{T} + 2\right)^d, \\
&\mathbb{E}_t = \left\{z \in \mathbb{Q}_t: \ |C_z^1|^{-1}|\{V_t = 1\} \cap C_z^1| \leq \frac{1}{2}\mu\right\},
\end{aligned} \tag{2.3}$$

label all blocks that intersect $Q_t$ and those where the "solid skeleton" covers a small fraction of volume.

If $\underline{\lim}_{t\to\infty} H_1(t) \geq H_*$, then for large $t$ and $z \notin \mathbb{E}_t$

$$\int_{C_z^1} |\phi|^2 \leq c\mu^{-1}\hat{H}^2 \int_{C_z^1} (|\nabla\phi|^2 + \beta\|V_t^{1/2}\phi\|_2^2), \quad \hat{H}^2 = \max\left\{H_*^2, \frac{1}{\beta}\right\}. \tag{2.4}$$

This estimate follows from inequality (1.17) with $Q^+ = C_z^1$, $Q^0 = C_z^1 \cap \{V_t = 1\}$, $Q^1 = Q^+ \setminus Q^0$, and $\alpha_*^d = 2/\mu - 1$. When $\mathbb{E}_t = \emptyset$ and $t$ is large, inequality (2.4) provides for PE (2.2) the rough lower bound

$$\lambda(t, \alpha, \beta) \geq c\mu\min\{H_*^{-2}, \beta\} \geq c_1\mu\min\{\mu^{2/d}, \beta\}.$$

**Lemma 2.1.** *If $H_1(t) \to H_* > (32\mu^{-2}d)^{1/d}$ as $t \to \infty$, then $\mathbf{P}\{\mathbb{E}_t = \emptyset\} \to 1$.*

**Proof.** By condition (1.5) $|\{V_t = 1\} \cap C_z^1| \geq \Xi_z$, where $\Xi_z = \sum_{C_\zeta^0 \subset C_z^1} \xi_\zeta$. Hence $A_z = \{|\{V_t = 1\} \cap C_z^1|/|C_z^1| \leq \frac{1}{2}\mu\} \subseteq \{\Xi_z \leq \frac{1}{2}\mu\}$.

For large $t$, it follows from (2.4), the inequality of S. N. Bernstein (see [5], Chapter 3.4), the restrictions on $H_1$, and the estimate $\mu = \mathbf{E}\xi_z \leq 1$, that $T = H_1/\tau > H_*/\sqrt{2}$ and

$$\mathbf{P}(A_z) \leq \mathbf{P}\left\{\Xi_z \leq \frac{1}{2}\mu\right\} \leq \exp\left\{-\frac{\mu^2 T^d/8}{1 + \mu/2}\right\} \leq \exp\left\{-\frac{1}{16}\mu^2 H_*^d \ln t\right\}.$$

Since $\ln\#(\mathbb{Q}_t) = (1 + \mathrm{o}(1))d\ln t$ by (2.3), this implies the relations

$$\overline{\lim_{t\to\infty}} \mathbf{P}\{\mathbb{E}_t \neq \emptyset\} = \overline{\lim_{t\to\infty}} \mathbf{P}\left(\bigcup_{z\in\mathbb{Q}_t} A_z\right) \leq \overline{\lim_{t\to\infty}} \#(\mathbb{Q}_t)\exp\left\{-\frac{1}{16}\mu^2 H_*^d \ln t\right\} = 0. \qquad \square$$

Below the size of blocks $H_1(t)$ satisfies the condition of Lemma 2.1, and $L = L(t)$ is a large natural number. The normalized PE (2.2) is estimated using its counterparts

$$\lambda_z^{(L)} = \inf\{\|\phi\|_{L^2}^{-2}(\mathcal{K}_\alpha(\phi) + \beta\|V_t^{1/2}\phi\|_2^2): \phi \in \mathbf{H}_0^1(Q_{z-}^{(L)})\} \tag{2.5}$$

for the same operator restricted to functions that vanish outside the cubes $Q_{z-}^{(L)} = LH_1(z + \frac{14}{10}Q)$ contained in larger cubes $Q_z^{(L)} = LH_1(z + \frac{15}{10}Q)$.

**Lemma 2.2.** *Define $\mathbb{J}_L = \{z \in \mathbb{Z}^d: Q_z^{(L)} \cap Q_t \neq \emptyset\}$. If $H_1$ satisfies the conditions of Lemma 2.1, then for each $\varepsilon > 0$ and $\hat{c} = \hat{c}(\alpha, \beta, H_*)$*

$$\underline{\lim_{t\to\infty}} \mathbf{P}\left\{\lambda(t, \alpha, \beta) \geq (1 - \varepsilon)\left(1 + \frac{\hat{c}}{L^2}\right)^{-1}\min_{z\in\mathbb{J}_t}\lambda_z^{(L)}\right\} = 1.$$

In Lemma 2.2 the ratio $\hat{c}/L^2$ can be made arbitrarily small by the choice of $L$. The random variables $\lambda_z^{(L)}$ are identically distributed.

**Proof of Lemma 2.2.** It suffices to derive the estimate of the lemma assuming that $\mathbb{E}_t = \emptyset$ because $\lim_{t\to\infty} \mathbf{P}\{\mathbb{E}_t = \emptyset\} = 1$.

Following [4], choose a smooth function $\zeta(x) \in [0, 1]$ such that $\zeta(x) = 1$ for $x \in Q$ and $\zeta(x) = 0$ for $x \notin \frac{14}{10}Q$, while $\sum_{z\in\mathbb{Z}^d}\zeta^2(x - z) = 1$ and $|\nabla\zeta(x)| \leq c$ for all $x \in \mathbb{R}^d$.

Define $\zeta_z(x) = \zeta((LH_1)^{-1}x - z)$ and take a function $\psi(x) \in \mathbf{H}_0^1(Q_t)$. By the above, the functions $\psi_z(x) = \zeta_z(x)\psi(x)$ satisfy the equalities (notation is that of (1.11))

$$\sum_{z\in\mathbb{Z}^d} |\psi_z(x)|^2 = |\psi(x)|^2,$$

$$\sum_{z \in \mathbb{Z}^d} \zeta_z^2(x) \widehat{K}_\alpha(\phi(x)) = \widehat{K}_\alpha(\phi(x)),$$

so simple calculations show that

$$\sum_{z \in \mathbb{Z}^d} \mathcal{K}_\alpha(\psi_z) = \mathcal{K}_\alpha(\psi) + \omega, \quad |\omega| \leq \int_{Q_t} |\psi|^2 \left(1 + \left(\frac{1}{\alpha}\right)\right) \sum_{z \in \mathbb{Z}^d} |\nabla \zeta_z|^2,$$

where $\sum_{z \in \mathbb{J}_t} |\nabla \zeta_z(x)|^2 \leq c(LH_1)^{-2}$ because each point of $\mathbb{R}^d$ belongs to a uniformly bounded number of sets $\{\nabla \zeta_z \neq 0\}$.

It follows from (2.5) that $\lambda_z^{(L)} \|\psi_z\|_2^2 \leq \mathcal{K}_\alpha(\psi_z) + \beta \int_{Q_z^{(L)}} V_t |\psi_z|^2$ for each $z$. Combined with estimate (2.4), the inequality shows that for each $\psi \in \mathbf{H}_0^1(Q_t)$

$$\min_{z \in \mathbb{J}_t} \lambda_z^{(L)} \|\psi\|_2^2 = \min_{z \in \mathbb{J}_t} \lambda_z^{(L)} \sum_{z \in \mathbb{J}_t} \|\psi_z\|_2^2$$

$$\leq \sum_{z \in \mathbb{J}_t} \left( \mathcal{K}_\alpha(\psi_z) + \beta \int_{Q_z^{(L)}} V_t |\psi_z|^2 \right) \leq \mathcal{K}_\alpha(\psi) + \beta \int_{Q_t} V_t |\psi|^2 + c(LH_1)^{-2} \|\psi\|_2^2$$

$$\leq (1 + c\hat{H}^2 (LH_1)^{-2}) \left( \mathcal{K}_\alpha(\psi) + \beta \int_{Q_t} V_t |\psi|^2 \right). \qquad \Box$$

*2.2. Reduction to individual Rayleigh quotients*

For the functions that determine the "partial PE" (2.5) used in Lemma 2.2, the norm $\|\nabla \phi\|_2$ is bounded from above and below by quantities that depend only on $L$, $\beta$, and $H_1$.

**Lemma 2.3.** *For large $t$ and $z \in \mathbb{J}_t$*

$$\lambda_z^{(L)} = \inf\{(\mathcal{K}_\alpha(\phi) + \beta \|V_t^{1/2} \phi\|_2^2) \colon \phi \in \Phi_z^{(L)}\}, \tag{2.6}$$

*where $\Phi_z^{(L)} = \{\phi \in \mathbf{H}_0^1(Q_{z-}^{(L)}) \colon \lambda_- \leq \|\nabla \phi\|_2^2 \leq 2\lambda_+, \|\phi\|_2^2 = 1\}$, $\lambda_- = \inf\{\|\phi\|_2^{-2} \|\nabla \phi\|_2^2 \colon \phi \in \mathbf{H}_0^1(Q_{z-}^{(L)})\}$, and $\lambda_+ = \beta + \inf\{\|\phi\|_2^{-2} \|\nabla \phi\|_2^2 \colon \phi \in \mathbf{H}_0^1(Q_{z-}^{(L)}), \operatorname{div}(\phi) = 0\}$.*

**Proof.** The inequality $\lambda_z^{(L)} \leq \lambda_+$ is immediate from definition (2.5) because $0 \leq V_t \leq 1$ and hence $\beta \|V_t^{1/2} \phi\|_2^2 \leq \beta \|\phi\|_2^2$. Consequently, if $\|\phi\|_2 = 1$ and $\|\nabla \phi\|_2^2 > 2\lambda_+$, then

$$\mathcal{K}_\alpha(\phi) + \beta \|V_t^{1/2} \phi\|_2^2 \geq \|\nabla \phi\|_2^2 \geq 2\lambda_+ > \lambda_z^{(L)}. \qquad \Box$$

**Lemma 2.4.** *If $H_1(t)$ satisfies the conditions of Lemma 2.1, then for each $\delta > 0$ and $t > t_*(\delta)$ there exists a finite family of test functions $G_\delta = G_\delta(H_*, L, \alpha, \beta, t) \subset \mathbf{H}_0^1(Q_0^{(L)})$ such that*

$$\lambda_0^{(L)} \geq \min_{\phi \in G_\delta} \{\|\phi\|_2^{-2} (\mathcal{K}_\alpha(\phi_\delta) + \beta \|V_t^{1/2} \phi\|_2^2)\} - c\delta.$$

*The value of the constant $c$ is determined by $H_*$, $L$, $\alpha$ and $\beta$. Neither the number of functions $\#(G_\delta)$ nor $c$ depends on $t$.*

**Proof.** The change of scale $g_*(x) = g(H_* H_1^{-1} x)$ reduces functions on $Q_0^{(L)} = \frac{15}{10} LH_1 Q$ to ones defined on the set $H_* H_1^{-1} Q_0^{(L)} = \frac{15}{10} LH_* Q$ that is independent of $t$. It is easily seen that $|\lambda_0^{(L)}/\lambda_* - 1| \leq c_1 |H_1 - H_*|$ for large $t$, where $\lambda_* = \inf_{\phi_* \in \Phi_0^{(L)}} \|\phi_*\|_2^{-2} (\mathcal{K}_\alpha(\nabla \phi_*) + \|V_*^{1/2} \phi_*\|_2^2)$ and the constant $c_1(H_*, L, \beta, \alpha)$ is independent of $t$.

Let $n(x) \geq 0$ be a $C^\infty$-smooth kernel such that $\int n(x)\,dx = 1$ and $n(x) = 0$ for $|x| \geq 1$. For functions $\phi \in \Phi_0^{(L)}$ and small $\delta > 0$, the convolutions $\phi_\delta(x) = \int \phi_*(x + \delta y)n(y)\,dy$ vanish outside $H_*H_1^{-1}Q_0^{(L)}$ by (2.5) and (2.6). It follows from well-known estimates for convolutions that $\|\phi_\delta\|_2 \leq \|\phi_*\|_2$ and

$$\|\phi_* - \phi_\delta\|_2 \leq c\lambda_+^{1/2}\delta\|\phi_*\|_2, \qquad \mathcal{K}_\alpha(\phi_\delta) \leq \mathcal{K}_\alpha(\phi_*).$$

Since $0 \leq V_t \leq 1$, the above formula implies that for each $\phi \in \Phi_0^{(L)}$

$$1 - c_1\delta \leq \|\phi_\delta\|_2^2 \leq 1, \qquad |\,\|V_*^{1/2}\phi_*\|_2^2 - \|V_*^{1/2}\phi_\delta\|_2^2\,| \leq c_2\delta,$$

where $c_i$ are positive constants. It is immediate that for small $\delta > 0$

$$\lambda_0^{(L)} \geq \inf_{\phi \in \Phi_0^{(L)}}(\mathcal{K}_\alpha(\phi) + \beta\|V_t^{1/2}\phi\|_2^2) \geq \inf_{\phi \in \Phi_0^{(L)}} \frac{\mathcal{K}_\alpha(\phi_\delta) + \beta\|V_*^{1/2}\phi_\delta\|_2^2}{\|\phi_\delta\|_2^2} - c_3\delta.$$

For each fixed $\delta > 0$, the function $\phi_\delta$ and all its derivatives are uniformly bounded by quantities proportional to $\|\phi\|_2$. Hence the unit ball in $L^2(Q_{0-}^{(L)})$, which contains $\Phi_0$, is transformed by change of scale and convolution into a pre-compact subset of $\mathbf{H}_0^1(\frac{15}{10}LH_*Q)$ containing $\Phi_0^{(L,\delta)} = \{\phi_\delta\colon \phi \in \Phi_0^{(L)}\}$. This ensures the existence of the set $G_\delta$ and the constant $c$ of the lemma. $\qquad\square$

**Lemma 2.5.** *If $H_1(t)$ satisfies the conditions of Lemma 2.1, then there exist positive constants $c_i$ such that for each fixed natural $L$ and $\delta > 0$ the scaled PE (2.2) satisfies the relation*

$$\forall w > 0 \quad \overline{\lim}_{t \to \infty} \mathbf{P}\{\lambda(t, \alpha, \beta) < w\} \leq \overline{\lim}_{t \to \infty} \#(\mathbb{J}_t)\#(G_\delta)\mathbf{q}^*(t),$$

*where the set $G_\delta$ and its cardinality $\#(G_\delta)$ are described in Lemma 2.4, and $\mathbf{q}^*(t) = \sup\{\mathbf{P}\{\|\phi\|_2^{-2}(\mathcal{K}_\alpha(\phi) + \beta\|V_t^{1/2}\phi\|_2^2) < w^0\}\colon \phi \not\equiv 0, \phi \in \mathbf{H}_0^1(\mathbb{R}^d)\}$ with $w^0 = (w + c_1\delta)(1 + \hat{c}/L^2).$*

**Proof.** If $\mathbb{E}_t = \emptyset$, then the estimate of Lemma 2.2 holds true, and $\mathbf{p} = \overline{\lim}_{t\to\infty} \mathbf{P}\{\lambda(t, \alpha, \beta) < w\}$ satisfies the inequality

$$\mathbf{p} \leq \overline{\lim}_{t \to \infty} \mathbf{P}\{\mathbb{E}_t \neq \emptyset\} + \overline{\lim}_{t \to \infty} \mathbf{P}\Big\{\min_{z \in \mathbb{J}_t} \lambda_z^{(L)} < w^*\Big\} \leq \overline{\lim}_{t \to \infty} \mathbf{P}\{\mathbb{E}_t \neq \emptyset\} + \overline{\lim}_{t \to \infty} \#(\mathbb{J}_t)\mathbf{P}\{\lambda_0^{(L)} < w^*\}$$

with $w^* = w(1 + \hat{c}/L^2)$ because all random variables $\lambda_z^{(L)}$ have the same distribution. The desired estimate follows from Lemma 2.4:

$$\mathbf{P}\{\lambda_0^{(L)} < w^*\} \leq \mathbf{P}\Big\{\min_{\phi \in G_\delta}\{\|\phi_\delta\|_2^{-2}(\mathcal{K}_\alpha(\phi_\delta) + \beta\|V_*^{1/2}\phi_\delta\|_2^2)\} < w^* + c\delta\Big\}$$

$$\leq \#(G_\delta)\sup_{\phi \in \mathbf{H}^1(\mathbb{R}^d)}\mathbf{P}\{\|\phi_\delta\|_2^{-2}(\mathcal{K}_\alpha(\phi_\delta) + \beta\|V_*^{1/2}\phi_\delta\|_2^2) < w^* + c\delta\}. \qquad\square$$

### 2.3. Feasibility of low individual Rayleigh quotients

*Rarity of small values of potential term*
The large deviation techniques used below are exposed in [5], Chapter 8.

Consider the averages $\langle f \rangle_z = \tau^{-d}\int_{C_z^0} f(x)\,dx$ over cells $C_z^0$ of size $\tau$. It is easy to see that $|\langle V_t|\phi|^2\rangle_z - \langle V_t\rangle_z\langle|\phi|^2\rangle_z| \geq -\langle\,|\phi|^2 - \langle|\phi|^2\rangle_z\rangle_z$ because $0 \leq V_t \leq 1$, while the Hölder inequality and the classical multiplicative inequalities for Sobolev spaces (see [3], Chapter II.2) imply the estimates

$$\langle|\,|\phi|^2 - \langle|\phi|^2\rangle_z|\rangle_z = \langle|\,|\phi|^2 - |\langle\phi\rangle_z|^2 - |\phi|^2 - \langle\phi\rangle_z|^2|\rangle_z,$$

$$\leq 2\langle|\,|\phi|^2 - |\langle\phi\rangle_z|^2|\rangle_z \leq c_1\tau\langle|\nabla\phi|^2\rangle_z^{1/2}\langle|\phi|^2\rangle_z^{1/2}. \tag{2.7}$$

By condition (1.5) $\xi_z \leq \langle V_t \rangle_z \leq 1$, which leads to the inequalities

$$\langle V_t |\phi|^2 \rangle_z \geq X_z - c_1 \tau \langle |\nabla \phi|^2 \rangle_z^{1/2} \langle |\phi|^2 \rangle_z^{1/2}, \quad 0 \leq X_z \stackrel{\text{def}}{=} \xi_z \langle |\phi|^2 \rangle_z \leq \langle |\phi|^2 \rangle_z,$$

where the random variables $X_z$ are independent. Summing them and applying the Cauchy inequality results in the estimate

$$\|V_t^{1/2} \phi\|_2^2 \geq \tau^d \Xi - c_1 \tau \|\nabla \phi\|_2 \|\phi\|_2, \quad \Xi \stackrel{\text{def}}{=} \sum_{z \in \mathbb{Z}^d} X_z. \tag{2.8}$$

For $\tilde{x}(t) = x(t) + c_1 \tau \|\nabla \phi\|_2 \|\phi\|_2$, this implies the inequality

$$\mathbf{P}\{\|V_t^{1/2} \phi\|_2^2 \leq x(t)\} \leq \mathbf{P}\{\tau^d \Xi \leq \tilde{x}(t)\}. \tag{2.9}$$

In notation of Theorem 1.2, the exponential moments of random variables of (2.8) are expressed through $G_z(h) = -\ln \mathbf{E} e^{-h X_z} = G(h \langle |\phi|^2 \rangle_z)$ and

$$-\ln \mathbf{E} \exp\{-h\Xi\} = \tau^{-d} \overline{G}_t(h; \phi),$$

$$\overline{G}_t(h; \phi) \stackrel{\text{def}}{=} \tau^d \sum_{z \in \mathbb{Z}^d} G_z(h) = \tau^d \sum_{z \in \mathbb{Z}^d} G(h \langle |\phi|^2 \rangle_z). \tag{2.10}$$

The functional $\overline{G}_t(h; \phi)$ is essentially a majorant for $-\ln \mathbf{E} \exp\{-h\|V_t^{1/2} \phi\|_2^2\}$.

Denote by $\widehat{X}_z$ independent random variables whose distributions are obtained from those of $X_z$ by the Cramér transform:

$$\mathbf{P}\{\widehat{X}_z \in B\} = F_{t,z,h}(B) = \frac{\mathbf{E} \exp\{-h X_z\} \mathbf{1}\{X_z \in B\}}{\mathbf{E} \exp\{-h X_z\}}, \quad h > 0. \tag{2.11}$$

The expectation of $\widehat{X}_z$ is $\mathbf{E} \widehat{X}_z = G_z'(h)$ (here and below prime stands for $\mathrm{d}/\mathrm{d}h$), and its variance satisfies the inequality

$$\mathbf{var}(\widehat{X}_z) = -G_z''(h) \leq \frac{\mathbf{E} X_z^2 \exp\{-h X_z\}}{\mathbf{E} \exp\{-h X_z\}} \leq \langle |\phi|^2 \rangle_z^2. \tag{2.12}$$

Set $\widehat{\Xi}_0 = \sum_{z \in \mathbb{Z}^d}(\widehat{X}_z - \mathbf{E} \widehat{X}_z)$. The standard inversion formula for the Cramér transform yields the equality

$$\mathbf{P}\{\Xi \leq \tau^{-d} x\} = \exp\{-\tau^{-d}(\overline{G}_t(h) - h\overline{G}_t'(h))\} \mathbf{E} e^{h\widehat{\Xi}_0} \mathbf{1}_A, \tag{2.13}$$

where $\mathbf{1}_A$ is the indicator of the event $A = \{\widehat{\Xi}_0 \leq \tau^{-d}(x - \overline{G}_t'(h))\}$.

**Lemma 2.6.** *Let $h > 0$ be a fixed positive number and $\phi \in \mathbf{H}^1(\mathbb{R}^d)$ a fixed function. If $\overline{\lim}_{t \to \infty}(x(t) - \overline{G}_t'(h)) < 0$, then*

$$\overline{\lim_{t \to \infty}} \, \tau^d \ln \mathbf{P}\{\|V_t^{1/2} \phi\|^2 \leq x(t)\} \leq -\varliminf_{t \to \infty} (\overline{G}_t(h) - h\overline{G}_t'(h)).$$

**Proof.** For large $t$, the assumptions of the lemma guarantee that $\tilde{x}(t) - \overline{G}'(h) < 0$ in (2.9), so by (2.13)

$$\ln \mathbf{P}\{\Xi \leq \tau^{-d} \tilde{x}(t)\} \leq -\tau^{-d}(\overline{G}_t(h) - h\overline{G}_t'(h)).$$

Indeed, the expectation in (2.13) does not exceed one because $\widehat{\Xi}_0 \leq 0$ over the domain of integration in (2.13). $\qquad \square$

The following technical lemma presents the estimate of Lemma 2.6 in a more tractable form (see the proof in Appendix A).

**Lemma 2.7.** *Let $h > 0$ be a fixed positive number and $\phi \in \mathbf{H}^1(\mathbb{R}^d)$ a fixed function. If $\overline{\lim}_{t\to\infty}(x(t) - \Gamma'(h;\phi)) < 0$, then*

$$\overline{\lim_{t\to\infty}} \tau^d \ln \mathbf{P}\{\|V_t^{1/2}\phi\|^2 \le x(t)\} \le -(\Gamma(h;\phi) - h\Gamma'(h;\phi)).$$

*Feasible values of individual Rayleigh quotients*

The typical values of the ratio $\|\phi\|_2^{-2}\|V_t^{1/2}\phi\|^2$ for a non-zero function $\phi \in \mathbf{H}^1(\mathbb{R}^d)$ are characterized by functional (1.15). To simplify calculations it is assumed that $\|\phi\|_2 = 1$ unless stated otherwise. Notation is that of (1.6).

The function $G(h)$ of (1.14) is non-negative, non-decreasing, concave, and bounded. By its definition $e^{-G(u)} \equiv \mathbf{E}\exp\{-u\xi_0\}$, and it follows from (1.5) that $p \le \mathbf{E}\exp\{-u\xi_0\} \le 1$ for all $u > 0$, so $0 = G(0) \le G(u) \le \nu$. Its derivatives can be expressed in terms of expectations (see (2.11) and (2.12)), and Lebesgue's theorem provides the limits as $h \to \infty$:

$$0 < G'(u) = e^{G(u)}\mathbf{E}\xi_0 e^{-u\xi_0} \le \frac{\mu}{p}, \qquad G'(0) = \mu,$$

$$\lim_{u\to\infty} G(u) = \nu, \qquad \lim_{u\to\infty} uG'(u) = 0. \tag{2.14}$$

Moreover, $-G'(u) \le G''(u) = (e^{G(u)}\mathbf{E}\xi_0 e^{-u\xi_0})^2 - e^{G(u)}\mathbf{E}\xi_0^2 e^{-u\xi_0} \le 0$ because $\xi_0^2 \le \xi_0$.

The function $\Gamma(h;\phi)$ of (1.14) is well defined and has derivatives

$$\Gamma'(h;\phi) = \int |\phi(x)|^2 G'(h|\phi(x)|^2)\,\mathrm{d}x > 0, \qquad \Gamma''(h) = \int |\phi|^4 G''(h|\phi|^2)\,\mathrm{d}x > 0,$$

for $\phi \in L^2(\mathbb{R}^d)$ because by the above $G$ has bounded derivative, so both $G(h|\phi|^2)$ and $|\phi|^2 G'(h|\phi|^2)$ are integrable. The second derivative exists for $h > 0$ because $|\phi|^4|G''(h|\phi|^2)| \le Ch^{-1}|\phi|^2 \in L^1(\mathbb{R}^d)$, $C = \sup_{u>0} uG'(u)$, by the estimate for $G''$ following (2.14). It is easily seen that

$$0 < \Gamma(h;\phi) < \left(\frac{\mu}{p}\right)h\|\phi\|_2^2, \qquad \lim_{h\to0+} \Gamma(h;\phi) = 0, \tag{2.15}$$

$$\lim_{h\to0+} \Gamma'(h;\phi) = \mu\|\phi\|_2^2, \qquad 0 < \Gamma'(h;\phi) \le \left(\frac{\mu}{p}\right)\|\phi\|_2^2. \tag{2.16}$$

If $\phi \in \mathbb{R}^d$ and $D_\infty = D_\infty(\phi) \stackrel{\text{def}}{=} \nu|\{|\phi| > 0\}| < \infty$, then

$$\Gamma(h;\phi) < \lim_{h\to\infty} \Gamma(h;\phi) = D_\infty, \qquad \lim_{h\to\infty} h\Gamma'(h;\phi) = 0. \tag{2.17}$$

**Lemma 2.8.** *The function $D \mapsto \mathcal{G}(\phi; D)$ of (1.15) is continuous and non-increasing on $[0,\infty)$. If $\phi$ vanishes outside a set of finite measure, then $\mathcal{G}(\phi; D) > 0$ for $D < D_\infty$ and $\mathcal{G}(\phi; D) = 0$ for $D \ge D_\infty$, where the threshold value $D_\infty$ is defined in (2.17).*

*For each $D \in (0, D_\infty]$ and arbitrarily small numbers $\varepsilon, \delta > 0$, one can find a finite number $h > 0$ such that*

$$0 \le \mathcal{G}(\phi; D) - h^{-1}(\Gamma(h) - D) \le \delta, \qquad D - \varepsilon \le \Gamma(h) - h\Gamma'(h) \le D. \tag{2.18}$$

*If $0 < D < D_\infty$, then there exists $h > 0$ such that (2.18) holds with $\delta = \varepsilon = 0$.*

**Proof.** If $D < D_\infty$, the function in (1.15) attains maximum at a single point. Indeed, its derivative can be represented in the form

$$(h^{-1}(\Gamma(h) - D))' = h^{-2}U(h), \qquad U(h) = h\Gamma'(h) - (\Gamma(h) - D).$$

Clearly, $U(0) = D > 0$ and $\lim_{h \to \infty} U(h) = D - D_\infty < 0$ (see (2.15) and (2.17)). For $h > 0$ the derivative $U'(h) = h\varGamma''(h)$ is strictly negative, so by the implicit function theorem the equation

$$U(h) \equiv D - (\varGamma(h) - h\varGamma'(h)) = 0$$

defines a unique strictly increasing function $h(D)$. This point is the required maximum and by the above

$$\mathcal{G}(\phi, D) = h^{-1}(\varGamma(h(D)) - D) = [h^{-1}(h\varGamma'(h) - U(h))]|_{h=h(D)} = \varGamma'(h(D)) > 0.$$

Relations (2.18) are obviously true for $h = h(D)$ with $\delta = \varepsilon = 0$.

By (2.16) and (2.17) $\varGamma(h) - h\varGamma'(h) < D_\infty$ for each finite $h > 0$, so it follows from (2.17) that $\lim_{D \nearrow D_\infty} h(D) = \infty$ and $\lim_{D \nearrow D_\infty} \mathcal{G}(\phi; D) = 0$.

If $D \geq D_\infty$, then it follows from (2.16) and (2.17) that $\mathcal{G}(\phi; D) = \lim_{h \to \infty} \frac{1}{h}(\varGamma(h; \phi) - D) = 0$ because $U(h) > 0$ for all finite $h$. Moreover, $\lim_{h \to \infty} U(h) = 0$ for $D = D_\infty$. Hence it suffices to take a large enough value of $h > 0$ to satisfy (2.18). $\qquad \square$

Consider $\varSigma_{\alpha,\beta}(D) \stackrel{\text{def}}{=} \inf\{\mathcal{K}_\alpha(\phi) + \beta\mathcal{G}(\phi; D) \colon \|\phi\|_2 = 1, \phi \in \mathbf{H}^1(\mathbb{R}^d)\}$ (cf. (1.11)). Later it will be essential that $\varSigma_{\alpha,\beta}(d) = \mathcal{C}_{\alpha,\beta}$ (see (1.13)).

**Lemma 2.9.** (a) *If $\|\phi\|_2 = 1$ and $0 \leq w < \mathcal{G}(\phi; D)$, then*

$$\overline{\lim_{t \to \infty}} (\ln t)^{-1} \ln \mathbf{P}\{\tau^{-d}\|V_t^{1/2}\phi\|_2^2 \leq w\} < -D.$$

(b) *If $\phi \in \mathbf{H}^1(\mathbb{R}^d)$ and $\|\phi\|_2 = 1$, then for each $\varepsilon \in (0, \varSigma_{\alpha,\beta}(D))$*

$$\overline{\lim_{t \to \infty}} (\ln t)^{-1} \ln \mathbf{P}\left\{\mathcal{K}_\alpha(\phi) + \beta\|V_t^{1/2}\phi\|_2^2 < \varSigma_{\alpha,\beta}(D) - \frac{\varepsilon}{\beta}\right\} < -D.$$

*The probabilities of the lemma are zero for $w < 0$ or $\varepsilon > \varSigma_{\alpha,\beta}(D)$.*

**Proof of Lemma 2.9.** By Lemma 2.8 there exists $h > 0$ such that $\mathcal{G}(\phi, D) = h^{-1}(\varGamma(h, \phi) - D)$ and $\varGamma(h, \phi) - h\varGamma'(h, \phi) = D$. Consequently, $w - \varGamma'(h) = w - \mathcal{G}(\phi, D) < 0$, so assertion (a) of the lemma follows from Lemma 2.7.

Assertion (b) follows from (a). Inequality $\varSigma_{\alpha,\beta}(D) \leq \mathcal{K}_\alpha(\phi) + \beta\mathcal{G}(\phi, D)$ follows from the definitions of $\varSigma_{\alpha,\beta}$ and $\mathcal{G}$, so in case of $\mathcal{G}(\phi, D) > 0$

$$\mathbf{P}\{\mathcal{K}_\alpha(\phi) + \beta\|V_t^{1/2}\phi\|_2^2 < \varSigma_{\alpha,\beta} - \varepsilon\} \leq \mathbf{P}\left\{\|V_t^{1/2}\phi\|_2^2 < \mathcal{G}(D) - \frac{\varepsilon}{\beta}\right\}. \qquad \square$$

*2.4. Proof of Theorem 1.2*

Choose a partition of $\mathbb{R}^d$ into blocks of size $H_1(t)$ satisfying the conditions of Lemma 2.1. Apply Lemmas 2.2 and 2.5 with $w = \mathcal{C}_{\alpha,\beta} - 2\varepsilon$ selecting $L$ and $\delta$ so that in notation of the latter lemma

$$\overline{\lim_{t \to \infty}} \mathbf{P}\{\lambda(t, \alpha, \beta) < w\} \leq \overline{\lim_{t \to \infty}} \#(\mathbb{J}_t)\#(G_\delta)\mathbf{q}^*(t),$$

where $\mathbf{q}^*(t) = \sup_\phi \mathbf{P}\{\|\phi\|_2^{-2}(\mathcal{K}_\alpha(\phi) + \beta\|V_t^{1/2}\phi\|_2^2) < w^0\}$ and for large $t$

$$w^0 = (w + c_1\delta)\left(1 + \frac{\hat{c}}{L^2}\right) < \mathcal{C}_{\alpha,\beta} - \varepsilon.$$

By the construction $\lim_{t \to \infty} (\ln t)^{-1} \ln(\#(\mathbb{J}_t)\#(G(\delta))) = d$. Since

$$\overline{\lim_{t \to \infty}} (\ln t)^{-1} \ln \mathbf{q}^*(t) < -d$$

by the estimate for $w^0$ and Lemma 2.9 (with $D = d$ and $\Sigma_{\alpha,\beta}(d) = \mathcal{C}_{\alpha,\beta}$), it follows that $\overline{\lim}_{t\to\infty} \mathbf{P}\{\lambda(t,\alpha, \beta) < w\} = 0$.

## 3. Low compressibility bound: case of large $\beta$

The bulk of this section is occupied by the proof of the following lemma. Theorem 1.1 is derived from it at the end of the section.

**Lemma 3.1.** *Consider* $\mathcal{C}_{\alpha,\beta}$ *defined in (1.13) and* $\mathcal{C}_\alpha$ *of (1.4). Under the conditions of Theorem 1.1* $\lim_{\beta\to\infty} \mathcal{C}_{\alpha,\beta} = (\nu/d)^{2/d}\mathcal{C}_\alpha$.

**Proof.** The argument below makes use of scaling properties of the analogues of PE's (1.4) for domains of arbitrary positive measure:

$$\inf\{\|\phi\|_2^{-2}\|\nabla\phi\|_2^2\colon \operatorname{div}(\phi) = 0, |\{|\phi| > 0\}| = u\} = u^{-2/d}\mathcal{S},$$

$$\inf\{\|\phi\|_2^{-2}\mathcal{K}_\alpha(\phi)\colon |\{|\phi| > 0\}| = u\} = u^{-2/d}\mathcal{C}_\alpha. \tag{3.1}$$

By Lemma 2.5 $\mathcal{G}(\phi; d) = 0$ if $D_\infty(\phi) = \nu|\{\phi \neq 0\}| \leq d$. For this reason definition (1.13) and (3.1) imply the inequality

$$\mathcal{C}_{\beta,\alpha} \leq s_* \stackrel{\text{def}}{=} \inf\left\{\|\phi\|_2^{-2}\mathcal{K}_\alpha(\phi)\colon |\{\phi \neq 0\}| \leq \frac{d}{\nu}\right\} = \left(\frac{\nu}{d}\right)^{2/d}\mathcal{C}_\alpha. \tag{3.2}$$

Thus, the test functions that determine the value of $\mathcal{C}_{\beta,\alpha}$ are, for each $\beta > 0$, in the set

$$\Phi_0 = \{\phi \in \mathbf{H}^1(\mathbb{R}^d)\colon \|\phi\|_2 = 1, \mathcal{K}_\alpha(\phi) \leq s_*\}. \tag{3.3}$$

To prove Lemma 3.1, it suffices to show that the value of $\lim_{\beta\to\infty} \mathcal{C}_{\beta,\alpha}$ is determined by test functions satisfying condition $\mathcal{G}(\phi; d) = 0$. To do so, set (3.3) is divided, for each $\beta$, into a few subsets defined using a small number $\delta > 0$ and a function $\varepsilon(\beta) > 0$ such that $\varepsilon(\beta) \searrow 0$ and $\beta\varepsilon(\beta) \nearrow \infty$ as $\beta \to \infty$. Below $\phi_\varepsilon(x) = \min\{|\phi(x)|, \varepsilon^{1/2}(\beta)\}$ and

$$\Psi_1(\beta) = \{\phi\colon \|\phi_\varepsilon\|_2^2 \leq \varepsilon^{1/2}(\beta)\|\phi\|_2^2\} \cap \Phi_0, \tag{3.4}$$

$$\Psi_2(\beta) = \left\{\phi\colon |\{|\phi(x)|^2 > \phi_\varepsilon^2(x)\}| \leq \frac{d+\delta}{\nu}\right\} \cap \Phi_0. \tag{3.5}$$

The elimination of irrelevant test functions is based on properties of $G$ and $\Gamma$ summarized in (2.14)–(2.17).

(a) If $\phi \in \Phi_0 \setminus \Psi_1(\beta)$, then one gets a lower bound for $\mathcal{G}(\phi; d)$ choosing $h(\beta) = \varepsilon^{-1}(\beta)$ in (1.15). It is easily seen that $G(u) \geq C\mu u$ for $0 < u \leq 1$ in (2.14). Since $h(\beta)\phi_\varepsilon^2(x) \leq h(\beta)(\varepsilon^{1/2}(\beta))^2 \leq 1$, this leads to the inequalities $\Gamma(h(\beta); \phi) \geq \int G(h(\beta)\phi_\varepsilon^2(x))\,dx \geq C\mu h(\beta)\|\phi_\varepsilon\|_2^2 \geq C\mu h(\beta)\varepsilon^{1/2}(\beta)$.

It follows from (1.13), (1.15) and (3.2) that for large $\beta$ the exclusion of test functions from $\Phi^0 \setminus \Psi_1(\beta)$ does not influence the value of $\mathcal{C}_{\beta,\alpha}$ because of the uniform lower bound

$$\forall\phi \in \Phi^0 \setminus \Psi_1(\beta) \quad \beta\mathcal{G}(\phi; d) \geq \beta\varepsilon^{1/2}(\beta)(C\mu - \varepsilon^{-1/2}(\beta)d) \nearrow \infty.$$

(b) If $\phi \in \Psi_1(\beta) \setminus \Psi_2(\beta)$, then $|\phi|^2 > \phi_\varepsilon^2 \geq \varepsilon(\beta)$ on a set of large measure. Since $\lim_{u\to\infty} G(u) = \nu$, one can choose $u$ so large that $G(u) > (1 - (\delta/d)^2)\nu$ and take $h(\beta) = u\varepsilon^{-1}(\beta)$ to see that

$$\Gamma(h(\beta)) \geq G(u)|\{|\phi|^2 > \phi_\varepsilon^2\}| \geq \left(1 - \left(\frac{\delta}{d}\right)^2\right)\frac{\nu(d+\delta)}{\nu} \geq d + \frac{1}{2}\delta.$$

Thus, also test functions from $\Psi_1 \setminus \Psi_2$ can be excluded when the infimum in (1.13) is calculated for large $\beta$ because for them

$$\beta \mathcal{G}(\phi; d) \geq \beta h^{-1}(\beta)(\Gamma(h; \phi) - d) \geq \frac{1}{2}\left(\frac{\delta}{u}\right) \beta \varepsilon(\beta) \to \infty.$$

(c) By (3.2) and (3.3) it suffices to show that there exist functions $\varkappa_*(\delta) > 0$ and $B(\kappa)$ such that $\lim_{\delta \to 0+} \varkappa_*(\delta) = 0$ and for $\beta > B(\varkappa_*(\delta))$

$$\{\phi \colon \|\phi\|_2^{-2} \mathcal{K}_\alpha(\phi) \geq s_* - \varkappa_*(\delta)\} \cap \Psi^0 \neq \emptyset, \quad \Psi^0 \stackrel{\text{def}}{=} \Psi_1(\beta) \cap \Psi_2(\beta). \tag{3.6}$$

In the proof of (3.6) the Rayleigh quotient of a function $\phi \in \Psi^0(\beta)$ is compared with that of $\zeta_U \phi$. The special cutoff function $\zeta_U$, constructed individually for each test function $\phi$, vanishes outside the set $U$ whose measure does not significantly exceed $d/\nu$, so $\|\zeta_U \phi\|_2^{-2} \|\nabla(\zeta_U \phi)\|_2^2$ cannot be much smaller than $s_*$ of (3.2) and (3.3).

The construction[2] of $U$ and $\zeta_U$ is described in detail in Appendix B. For a test function $\phi \in \Psi^0$ it is based on the set

$$E = \{x \colon |\phi(x)| > \phi_\varepsilon(x)\}. \tag{3.7}$$

The set $U = U(\phi)$ (see (B.4)) results from approximation of $E$ by a finite union of cubes $C_z^k = H_k(z + Q)$ of sizes $H_k = T^{k-K_*} h_*$, $k = 0, \ldots, K_* - 1$, that are "empty;" i.e., $|C_z^k \cap E| > (1 - \gamma)|C_z^k|$. The cutoff function $\zeta_U$ is described in (B.7).

The parameters used in the construction of $U$ and $\zeta_U$ are functions of the number $\delta > 0$ in definition (3.5) such that as $\delta \to 0$

$$T = T(\delta) \nearrow \infty, \qquad \gamma = T^{-\kappa}, \qquad m_0 \sim T^\kappa \in \mathbb{N}, \qquad 0 < h_* < \frac{1}{2}\gamma, \tag{3.8}$$

where $T(\delta)$ is a large odd natural number and $\kappa > 0$ a fixed small positive exponent. The natural number $K_* = K_*(T)$ is chosen so as to satisfy the conditions $m_0/T < \gamma$ and $m_0^d/(K_* + 1) < \gamma$ (see (B.5)).

Since $\|\phi\|_2 = 1$, it follows from (B.10) and (B.6) that for $\varepsilon_1 = (\alpha \gamma m_0^2)^{-1/2}$

$$(1 + c_1 \varepsilon_1)\left(\frac{|U|^{2/d}}{\mathcal{C}_\alpha}\right) \int_U K_\alpha(\phi(x)) \, \mathrm{d}x \geq 1 - c_3(\varepsilon_1 H_{k_0}^{-2} + \gamma^{-1}) \|\phi_\varepsilon\|_2^2 - c_4 \gamma \|\nabla \phi\|_2^2$$

and $|U| \leq (1 + c_4 \gamma)|E| \leq (1 + c_4 \gamma)(d + \delta)/\nu$. Combined, these inequalities yield the estimate

$$\|\phi\|_2^{-2} \mathcal{K}_\alpha(\phi) \geq \mathcal{C}_\alpha |U|^{-2/d}(1 - c_5(\varepsilon_1 + \gamma) - c_6(\varepsilon_1 H_{k_0}^{-2} + \gamma^{-1}) \|\phi_\varepsilon\|_2^2)$$

$$\geq \mathcal{C}_\alpha \left(\frac{\nu}{d}\right)^{2/d} \left(1 + \frac{\delta}{d}\right)^{-2/d} (1 - c'(\varepsilon_1 + \gamma) - c''(\varepsilon_1 H_{k_0}^{-2} + \gamma^{-1}) \|\phi_\varepsilon\|_2^2). \tag{3.9}$$

In the latter estimate, $\|\phi_\varepsilon\|_2^2 \leq \varepsilon \nu^{-1}(d + \delta)$ converges to zero uniformly on $\Psi^0$ as $\beta \to \infty$, while the quantities $\gamma^{-1} H_{k_0}^2$ and $\varepsilon_1$, which do not depend on $\beta$, can be made arbitrarily small if $T = T(\delta)$ is chosen large enough.

Thus, given $\delta$, one can select $T = T(\delta)$ so that (3.9) guarantees the estimate

$$\|\phi\|_2^{-2} \mathcal{K}_\alpha(\phi) > \mathcal{C}_\alpha \left(\frac{\nu}{d}\right)^{2/d} \Big/ (1 + 2\delta)$$

if $\Psi^0(\beta) \neq \emptyset$ and $\beta$ is large enough.

The arguments in (a)–(c) show that $\underline{\lim}_{\beta \to \infty} \mathcal{C}_{\alpha, \beta} = \mathcal{C}_\alpha$. The lemma is proved. $\qquad \square$

---

[2] It is, in essence, a simpler version of one used in [10] (and goes back to "the method of enlargement of obstacles" of [6, 7]).

**Proof of Theorem 1.1.** By Lemma 3.1 $\lim_{\beta \to 0+} \mathcal{C}_{\alpha,\beta} = \mathcal{C}_\alpha$. Combined with (1.11) and (1.4), this fact allows one to conclude that for each $\varepsilon > 0$ and sufficiently small values of $\alpha = \alpha(\varepsilon)$ and $\beta = \beta(\varepsilon)$.

$$\lim_{t \to \infty} \mathbf{P}\left\{ (\ln t)^{2/d} \mathfrak{S}_t > \left(\frac{\nu}{d}\right)^{2/d} \mathcal{S} - \varepsilon \right\} \geq \lim_{t \to \infty} \mathbf{P}\left\{ (\ln t)^{2/d} \mathfrak{C}_{\alpha,\beta,t} > \left(\frac{\nu}{d}\right)^{2/d} \mathcal{S} - \frac{\varepsilon}{2} \right\} = 1. \qquad \square$$

## Appendix A. Proof of Lemma 2.7

It suffices to establish the convergence of the derivative $\overline{G}'_t(h) \to \Gamma'(h)$ for each fixed value of $h > 0$; the convergence of $\overline{G}_t(h)$ follows. By definitions (1.14) and (2.10),

$$\overline{G}'_t(h) - \Gamma'(h) = \tau^d \sum_{z \in \mathbb{Z}^d} (G'_z(h) - \Gamma'_z(h)), \qquad \Gamma_z(h) = \langle G(h|\phi|^2) \rangle_z,$$

and it is easily seen that

$$G'_z(h) - \Gamma'_z(h) = \langle \mathbf{E}\xi_z U(\theta) W(\theta) \rangle_z \big|_{\theta=0}^{\theta=1}, \tag{A.1}$$

where $U(x,\theta) = ((1-\theta)|\phi(x)|^2 + \theta\langle|\phi|^2\rangle_z)$, and for each $x$ and $\theta$ the weight $W(x,\theta) = e^{-h\xi_z U(x,\theta)}/\mathbf{E}e^{-h\xi_z U(x,\theta)}$ satisfies the condition $\mathbf{E}W(\theta) = 1$. Moreover, $\dot{U}(x) \overset{\text{def}}{=} (\partial/\partial\theta)U = -(|\phi(x)|^2 - \langle|\phi|^2\rangle_z)$.

(a) For $d = 2, 3$ it suffices to use Barrow's formula in the form

$$G'_z(h) - \Gamma'_z(h) = \int_0^1 \left(\frac{\mathrm{d}}{\mathrm{d}\theta}\right) \langle \mathbf{E}\xi_z U(\theta) W(\theta) \rangle_z \, \mathrm{d}\theta = \int_0^1 \sum_{k=1}^3 \Psi_k \, \mathrm{d}\theta, \tag{A.2}$$

where $\Psi_1(\theta) = \langle \mathbf{E}\xi_z \dot{U}(\theta) W(\theta) \rangle_z$, $\Psi_2(\theta) = -h\langle \mathbf{E}\xi_z^2 U(\theta) \dot{U}(\theta) W(\theta) \rangle_z$, and $\Psi_3(\theta) = -h\langle (\mathbf{E}\xi_z \dot{U}(\theta) W(\theta))(\mathbf{E}\xi_z U(\theta) \times W(\theta)) \rangle_z$.

Indeed, both $U \geq 0$ and $|\dot{U}|$ can be estimated from above by non-random quantities, and the same estimates hold true for the cell averages of the expectations with weight $W$ in (A.2). It follows from (2.7), (1.16), and the Hölder inequality that

$$|\langle \mathbf{E}\xi_z U(\theta) W(\theta) \rangle_z| \leq \langle (|\phi|^2 + \langle|\phi|^2\rangle_z) \rangle_z \mathbf{E}\xi_z W(\theta) \leq 2\langle|\phi|^2\rangle_z,$$

$$\langle \mathbf{E}|\xi_z \dot{U}(\theta)| W(\theta) \rangle_z \leq \langle ||\phi|^2 - \langle|\phi|^2\rangle_z| \mathbf{E}\xi_z W(\theta) \rangle_z \leq C\tau \langle|\nabla\phi|^2\rangle_z^{1/2} \langle|\phi|^2\rangle_z^{1/2}. \tag{A.3}$$

A similar calculation shows that

$$|\langle (\mathbf{E}\xi_z \dot{U}(\theta) W(\theta))(\mathbf{E}\xi_z U(\theta) W(\theta)) \rangle_z| \leq \langle |\langle|\phi|^2\rangle_z - |\phi|^2|(|\phi|^2 + \langle|\phi|^2\rangle_z) \rangle_z \leq c\langle|\phi - \langle\phi\rangle_z|^2\rangle_z^{1/2} \langle|\phi|^6\rangle_z^{1/2}.$$

Combining the above estimates, one arrives from (A.2) at the inequality

$$|G'_z(h) - \Gamma'_z(h)| \leq C\langle|\phi - \langle\phi\rangle_z|^2\rangle_z^{1/2}((1+h)\langle|\phi|^2\rangle_z^{1/2} + h\langle|\phi|^6\rangle_z^{1/2}).$$

Application of (1.16) and the Cauchy inequality to sums of cell averages results in the estimate

$$|\overline{G}'_t(h) - \Gamma'(h)| \leq c\tau\|\nabla\phi\|_2(\|\phi\|_2 + \|\phi\|_6^3) \leq c\tau(\|\phi\|_{\mathbf{H}^1}^2 + \|\phi\|_{\mathbf{H}^1}^4). \tag{A.4}$$

(b) For $d \geq 4$, there is no estimate for $\|\phi\|_6$ in terms of $\|\phi\|_{\mathbf{H}^2}$. This difficulty is circumvented by truncation: the $\mathbb{R}^d$-valued test function $\phi \in \mathbf{H}^1(\mathbb{R}^d)$ is approximated by the function $\phi_r = (\phi_r^{(j)})$ with coordinates

$$\phi_r^{(j)}(x) = \min\{|\phi^{(j)}(x)|, r\} \operatorname{sign}(\phi^{(j)}(x)).$$

Clearly, $|\phi_r| \leq |\phi|$ and $|\phi_r - \phi| \leq |\phi|$. If $\phi \in \mathbf{H}^1(\mathbb{R}^d)$, then $\phi_r$ belongs to the same space and (see [3], Chapter II.3 or [2], Chapter 7.4)

$$\|\phi_r\|_p \leq \|\phi\|_p, \quad 1 \leq p < \infty, \qquad \|\nabla \phi_r\|_2 \leq \|\nabla \phi\|_2. \tag{A.5}$$

Define $U_r(x, \theta) = \xi_z((1-\theta)|\phi_r(x)|^2 + \theta \langle |\phi_r|^2 \rangle_z)$. Using the same weights $W$, difference (A.1) is represented in the form

$$G_z'(h) - \Gamma_z'(h) = A(\theta)|_{\theta=0}^{\theta=1} + \varepsilon_{1,z}(\theta)|_{\theta=0}^{\theta=1}, \tag{A.6}$$

where $A(\theta) = \langle \mathbf{E}\xi_z U_r(\theta) W(\theta) \rangle_z$ and $\varepsilon_{1,z}(\theta) = \langle \mathbf{E}\xi_z (U(\theta) - U_r(\theta)) W(\theta) \rangle_z$.

($b_1$) Evidently, $0 \leq |\phi|^2 - |\phi_r|^2 \leq 2|\phi - \phi_r| \, |\phi|$, so $\xi_z |U(x, \theta) - U_r(x, \theta)| \leq 2|\phi - \phi_r| \, |\phi| + \langle 2|\phi - \phi_r| \, |\phi| \rangle_z$ by (1.5) and

$$|\varepsilon_{1,z}(\theta)| \leq |\langle \mathbf{E}|U(\theta) - U_r(\theta)|W(\theta) \rangle_z| \leq c \langle |\phi - \phi_r|^2 \rangle_z^{1/2} \langle |\phi|^2 \rangle_z^{1/2}. \tag{A.7}$$

Choose $p \in (2, 2/(1 - 2/d)]$. Note that $\phi(x) - \phi_r(x) \neq 0$ only if $|\phi(x)| > r$. It follows from (1.16) that

$$\|\phi - \phi_r\|_2^2 \leq r^{-(p-2)} \int |\phi(x)|^p \, \mathrm{d}x \leq c(p) r^{-(p-2)} \|\phi\|_{\mathbf{H}^1}^p.$$

The Cauchy inequality for sums and the above inequality show that for $p = 2/(1 - 2/d)$

$$\tau^d \sum_{z \in \mathbb{Z}^d} |\varepsilon_{1,z}(\theta)| \leq c \|\phi - \phi_r\|_2 \|\phi\|_2 \leq c(d) r^{-2/(d-2)} \|\phi\|_{\mathbf{H}^1}^{d/(d-2)} \|\phi\|_2.$$

($b_2$) To estimate the first term on the right-hand side in (A.6), it proves convenient to represent it in the form

$$A(\theta)|_{\theta=0}^{\theta=1} = \int_0^1 \left(\frac{\mathrm{d}}{\mathrm{d}\theta}\right) \langle \xi_z U_r(\theta) W(\theta) \rangle_z \, \mathrm{d}\theta = \int_0^1 (F_1(\theta) + F_2(\theta) + F_3(\theta)) \, \mathrm{d}\theta,$$

where $F_1(\theta) = \langle \mathbf{E}\xi_z \dot{U}_r(\theta) \xi_z W(\theta) \rangle_z$, $F_2(\theta) = -h \langle \mathbf{E}\xi_z^2 U_r(\theta) \dot{U}(\theta) W(\theta) \rangle_z$, and $F_3(\theta) = -h \langle (\mathbf{E}\xi_z \dot{U}(\theta) W(\theta)) \times (\mathbf{E}\xi_z U_r(\theta) W(\theta)) \rangle_z$.

The argument used to derive (A.3) and inequality (A.5) prove that

$$|F_1(\theta)| \leq c \langle |\phi_r - \langle \phi_r \rangle_z|^2 \rangle_z^{1/2} \langle |\phi_r|^2 \rangle_z^{1/2} \leq c\tau \langle |\nabla \phi|^2 \rangle_z^{1/2} \langle |\phi|^2 \rangle_z^{1/2}.$$

It is immediate from (A.1) that $|U_r(\theta)| \leq cr^2$, so by (A.3)

$$|F_2(\theta)| \leq cr^2 h \langle \mathbf{E}|\dot{U}(\theta)|W(\theta) \rangle_z \leq cr^2 h\tau \langle |\nabla \phi|^2 \rangle_z^{1/2} \langle |\phi|^2 \rangle_z^{1/2}.$$

Finally, $|F_3(\theta)| \leq h \langle |\dot{U}(\theta)| cr^2 \rangle_z \leq cr^2 h\tau \langle |\nabla \phi|^2 \rangle_z^{1/2} \langle |\phi|^2 \rangle_z^{1/2}$.

It follows that for each cell

$$|A(\theta)|_{\theta=0}^{\theta=1}| \leq c\tau (1 + h)(1 + r^2) \langle |\nabla \phi|^2 \rangle_z^{1/2} \langle |\phi|^2 \rangle_z^{1/2}.$$

Combining the above estimate with (A.7) and applying the Cauchy inequality to the sums in $z$, one concludes that

$$|\overline{G}_t'(h) - \Gamma'(h)| \leq \tau^d \sum_{z \in \mathbb{Q}_t} (|\varepsilon_{1,z}| + c\tau \langle |\nabla \phi|^2 \rangle_z^{1/2} (1 + hr^2 \langle |\phi|^2 \rangle_z^{1/2}))$$

$$\leq c(r^{-2/(d-2)} \|\phi\|_{\mathbf{H}^1}^{d/(d-2)} \|\phi\|_2 + \tau(1 + h)(1 + r^2) \|\nabla \phi\|_2 \|\phi\|_2).$$

The choice $r = \tau^{-1/(d-1)}$ yields the final estimate which proves the lemma,

$$|\overline{G}_t'(h) - \Gamma'(h)| \leq c\tau^{1/(d-1)}(1 + h)(\|\phi\|_{\mathbf{H}^1}^{2+2/(d-2)} + \|\phi\|_{\mathbf{H}^1}^2).$$

## Appendix B. Block approximation of sets

*Construction of block approximation*

Below, a set $E$ of finite Lebesgue measure $|E|$ is approximated by the union of disjoint cubic blocks of several standard sizes as in [11].

An odd natural number $T$ and a number $H_0 > 0$ determine a hierarchy of scales and partitions of the space into equal *cells* $C_z^0$ or *blocks* of levels $k > 0$

$$C_z^k = H_k(z + Q), \quad H_k = T^k H_0, z \in \mathbb{Z}^d, \quad k = 0, 1, 2, \ldots.$$

Each $(k+1)$-block $C_z^{k+1}$ is the union of $T^d$ disjoint sub-blocks of level $k$. The number $T$ will be large later on.

One more parameter $\gamma \in (0,1)$ is used to sort blocks according to the fraction of their volume covered by $E$. The "empty" blocks $C_z^k$ are numbered by integer vectors

$$z \in \mathbb{E}_k = \{z : |C_z^k|^{-1}|C_z^k \cap E| > 1 - \gamma\}. \tag{B.1}$$

Since $|E|$ is finite, $\mathbb{E}_k = \emptyset$ for sufficiently large $k$.

Sets (B.1) define the solid sets $E_k^+ = \bigcup_{z \in \mathbb{E}_k} C_z^k$, $E_k = E \cap E_k^+$, and

$$\Psi_k = \bigcup_{\ell \geq k} E_\ell^+, \qquad \Phi_k = E_k^+ \setminus \Psi_{k+1}, \quad k = 0, 1, 2, \ldots. \tag{B.2}$$

Obviously, $\bigcup_{k \geq 0} E_k \subseteq E$ and $|E_k| \leq |E_k^+| \leq (1-\gamma)^{-1}|E_k|$.

By the construction, the sets $\Phi_k$ are disjoint, each set $\Phi_k$ is the union of empty blocks of level $k$, and $\Psi_k$ is the union of disjoint empty blocks $C_\zeta^\ell$, $\ell \geq k$, so $|\Phi_k| \leq (1-\gamma)^{-1}|\Phi_k \cap E|$ and $|\Psi_k| \leq (1-\gamma)^{-1}|E|$.

For each $k$, the complement of $\Psi_k$ is covered by non-empty blocks of level $k$ because all empty cubes of this and higher levels are included in $\Phi_l$, $l \geq k$.

To define the set $U$, choose one more parameter $m_0 \in \mathbb{N}$ so that $m_0 < T$.

Suppose that $\mathbb{E}_{K^*} \neq \emptyset$ for some $K^* > 0$. Let $k_0 \in \{0, 1, \ldots, K^*\}$ be the highest level such that

$$|\Phi_{k_0} \cap E| \leq \frac{|E|}{K^* + 1}, \qquad |\Phi_{k_0}| \leq \frac{(1-\gamma)^{-1}|E|}{K^* + 1}. \tag{B.3}$$

(Its existence follows from the inequality $\sum_{k=0}^{K^*} |\Phi_k \cap E| \leq |E|$.) The set $U$ is the $m_0 H_{k_0}$-neighborhood of the set $\Psi_{k_0}$ defined in (B.2):

$$U = \{x : \text{Dist}(x, \Psi_{k_0}) \leq m_0 H_{k_0}\}. \tag{B.4}$$

Its measure admits a simple estimate in terms of the original set $E$. If a point of $U$ is no further than $m_0 H_{k_0}$ from some empty block $C_\zeta^\ell \subset \Psi_{k_0}$, $\ell > k$, of size $H_\ell = T^{\ell-k_0} H_{k_0}$, it belongs to the larger cube $(1 + cm_0/T)^d H_\ell(\zeta + Q)$ whose volume is $(1 + cm_0/T)^d |C_\zeta^\ell|$. The volume of the set $\Phi_{k_0}$ is small by (B.3), and its $m_0 H_{k_0}$-neighborhood cannot have volume greater than $(2m_0 + 1)^d |\Phi_{k_0}|$. Consequently,

$$|U| \leq \left(1 + \frac{c'm_0}{T}\right)|\Psi_{k_0}| + c''m_0^d|E|(K^* + 1)^{-1}(1-\gamma)^{-1}.$$

This inequality is applied in situations where the parameters satisfy the condition

$$\frac{m_0}{T} < \gamma, \qquad \frac{m_0^d}{K_* + 1} < \gamma. \tag{B.5}$$

In this case, there exist positive constants $c_i$ such that for $\gamma < c_0$

$$|U| \leq \left(1 + c_1\left(\frac{m_0}{T} + \frac{m_0^d}{K^* + 1}\right)\right)\frac{|E|}{1 - \gamma} \leq (1 + c_2\gamma)|E|. \tag{B.6}$$

*Cutoff for block approximation*

Let $U$ be the set defined in (B.4) for a given value of $\gamma > 0$ (see (B.1)). Take a cutoff function $\zeta_U : \mathbb{R}^d \to [0, 1]$ such that

$$\zeta_U(x) = \begin{cases} 1, & x \in \Psi_{k_0}, \\ 0, & \text{Dist}(x, \Psi_{k_0}) \geq m_0 H_{k_0}, \end{cases} \qquad |\nabla \zeta_U(x)| \leq \frac{c}{m_0 H_{k_0}}, \tag{B.7}$$

where the constant $c$ does not depend on the shape of the set $\Psi_{k_0}$. All blocks $C_z^{k_0}$ where the gradient $\nabla \zeta_U$ does not vanish identically are non-empty.

It is proved in [11] that there exist constants $c_i > 0$ (independent of $\alpha > 0$), such that for $\widehat{K}_\alpha$ of (1.11), $\gamma m_0^2 > c_0$, and $\varepsilon_1 \stackrel{\text{def}}{=} (\alpha\gamma m_0^2)^{-1/2}$

$$\int_{\Psi_{k_0}} |\psi|^2 \leq (1 + c_1\varepsilon_1)\left(\frac{|U|^{2/d}}{\mathcal{C}_\alpha}\right)\int_U \widehat{K}_\alpha(\psi) + c_2\varepsilon_1 H_{k_0}^{-2}\|\psi\|_{L^2(\mathbb{R}^d\setminus E)}^2. \tag{B.8}$$

The proof is based on the fact that in each non-empty block $C_z^{k_0}$ the subsets $Q^0 = C_z^{k_0} \setminus E$ and $Q^1 = C_z^{k_0} \cap E$ have comparable measures: $|Q^1|/|Q^0| \leq 1/\gamma$. Hence, the Poincaré–Friedrichs inequality (1.17) can be used to show that

$$\|\phi\|_{L^2(C_z^{k_0})}^2 \leq c\gamma^{-1}\left(\int_{C_z^{k_0}\setminus E} |\phi|^2 + H_{k_0}^2\int_{C_z^{k_0}} |\nabla\phi(x)|^2\,\mathrm{d}x\right). \tag{B.9}$$

The inequality remains true for unions of disjoint non-empty $k_0$-blocks.

**Remark B.1.** *The proof of Lemma 3.1 uses inequality* (B.8) *in the form*

$$(1 + c_1\varepsilon_1)\left(\frac{|U|^{2/d}}{\mathcal{C}_\alpha}\right)\int_U K_\alpha(\phi(x))\,\mathrm{d}x \geq \|\phi\|_2^2 - c_3(\varepsilon_1 H_{k_0}^{-2} + \gamma^{-1})\|\phi_\varepsilon\|_2^2 - c_4\gamma\|\nabla\phi\|_2^2. \tag{B.10}$$

In this case, the selection of "empty" blocks is based on set (3.7). This set and parameters (3.8) of its block approximation are determined by the choice of a test function $\phi = \phi_\varepsilon$ with specified properties (see (3.4), (3.6), (3.9)).

The existence of empty blocks at level $K_*$ (of size $h_* < \gamma/2$) used in the construction is guaranteed by definition (3.3). Indeed, for a "non-empty" block $C_z^{K_*}$ the choice of $E$ and $h_*$ and inequality (B.9) provide the estimate $\int_{C_z^{K_*}} |\phi|^2 \leq c(\gamma^{-1}\int_{C_z^{K_*}} \phi_\varepsilon^2 + \gamma\int_{C_z^{K_*}} |\nabla\phi|^2)$. In absence of empty $K_*$-blocks this would imply the inequality $\|\phi\|_2^2 \leq c(\gamma^{-1}\|\phi_\varepsilon\|_2^2 + \gamma\|\nabla\phi\|_2^2)$ and the lower bound $\|\phi\|_2^{-2}\mathcal{K}_\alpha(\phi) \geq \gamma^{-1}\hat{c}(1 - \gamma^{-1}\varepsilon^{1/2}(\beta))$ for the Rayleigh quotient, which contradicts (3.3) if $\gamma = \gamma(T)$ is small and $\beta$ is large.

The special choice of set (3.7) and (B.9) used to define $U$ of (B.4) and $\zeta_U$ of (B.7) ensure the estimate $\int_{\mathbb{R}^d\setminus\Psi_{k_0}} |\phi|^2 \leq c(\gamma^{-1}\|\phi_\varepsilon\|_2^2 + \gamma\|\nabla\phi\|_2^2)$. (Note that $\gamma^{-1}H_k^2 \leq \gamma$ for $k \leq K_*$ by (3.8).) By inequality (B.8) with $\varepsilon_1 = (\alpha\gamma m_0^2)^{-1/2}$ and the choice of set (3.7) this leads to (B.10):

$$(1 + c_1\varepsilon_1)\left(\frac{|U|^{2/d}}{\mathcal{C}_\alpha}\right)\int_U K_\alpha(\phi) \geq \int_{\Psi_{k_0}} |\phi|^2 - c_2\varepsilon_1 H_{k_0}^{-2}\|\phi\|_{L^2(\mathbb{R}^d\setminus E)}^2$$

$$\geq \|\phi\|_2^2 - c_3(\varepsilon_1 H_{k_0}^{-2} + \gamma^{-1})\|\phi_\varepsilon\|_2^2 - c_4\gamma\|\nabla\phi\|_2^2.$$